
\input amstex
\documentstyle{amsppt}
\magnification=\magstep1
\baselineskip=12pt
\hsize=6truein
\vsize=8truein

\topmatter
\title On the rigidity for conformally compact Einstein manifolds
\endtitle

\author Jie Qing \\
Department of Mathematics \\
UCSC
\endauthor

\leftheadtext{Conformally compact Einstein}
\rightheadtext{Jie Qing}

\address Jie Qing, Dept. of Math., UC,
Santa Cruz, Santa Cruz, CA 95064. \endaddress

\email qing{\@}math.ucsc.edu \endemail

\abstract In this paper we prove that a conformally compact Einstein
manifold with the round sphere as its conformal infinity has to be the
hyperbolic space.
We do not assume the manifolds to be spin, but our 
approach relies on the positive mass theorem for asymptotic flat
manifolds.
The proof is based on understanding of positive eigenfunctions and
compactifications obtained by positive eigenfunctions. 
\endabstract

\endtopmatter
\document

\vskip 0.1in
\head 0. Introduction \endhead

In this paper we study the rigidity problem for conformally compact 
Einstein manifolds with the round sphere as their conformal infinity.
Quite recently there has been a great deal of interest in both physics
and mathematics community in the so-called Anti-de-Sitter/Conformal
Field Theory (in short AdS/CFT) correspondence. Conformally compact 
Einstein manifolds play an essential role in this correspondence.
In mathematics it
has been known for a long time that there are close connections
between the geometry of Minkowski space-time,
hyperbolic space and the round sphere.
Notably in a seminal
paper [FG] Fefferman and Graham
showed this approach to be very fruitful in conformal geometry. 

In establishing scalar curvature rigidity for asymptotically 
hyperbolic manifold as a generalization from the work of Witten [W]
on positive mass theorem for asymptotically
flat manifolds which are spin, in [AD], Andersson and Dahl proved that, 
if a conformally
compact Einstein manifold with the round sphere as its conformal infinity
is spin, then it has to be a hyperbolic space (please also see 
recent related works of X. Zhang [Z], Chru\'{s}ciel and Herzlich [CH],
and X. Wang [Wa]). It opens an interesting
question whether the spin structure is necessary to assure the rigidity.
There is some progress made by Anderson in [An1].

Before we state our results. let us briefly introduce what is a conformally
compact Einstein manifold. Let $X^{n+1}$ be a $n+1$ dimensional compact
manifold with boundary $M^n = \partial X$. $(X,g)$ is said to be a conformally
compact Einstein manifold if $\text{Ric}(g) = - n g$ and
$(X, s^2 g)$ is a compact Riemannian manifold with boundary, 
where $s$ is a defining function of the boundary $M$. 
Clearly the restriction of $s^2 g$ to $TM$ is a metric $\hat g$ on the 
boundary and $\hat g$ rescales upon changing the defining function $s$. 
Thus $(M, [\hat g])$ is determined by $(X, g)$ and called the conformal
infinity of $(X,g)$. 

\proclaim{Theorem 0.1} Suppose that $(X^{n+1},g)$ is a conformally compact 
Einstein manifold with the round sphere as its conformal infinity, and
$3\leq n\leq 6$. Then $(X,g)$ has to be the hyperbolic space.
\endproclaim

One simple yet very interesting calculation leading to Theorem 0.1
is the following. 

\proclaim{Lemma 0.2} Suppose that $(X, g)$ is a conformally compact
Einstein manifold. And suppose that $u$ is a positive eigenfunction,
i.e.  $\Delta u = (n+1) u$. Then $(X, u^{-2}g)$ is with scalar curvature
$$
R = n(n+1) (u^2 - |du|^2).
\tag 0.1
$$
\endproclaim

Here $\Delta$ is the trace of the Hessian in metric $g$.
Combining with the Bochner formula for eigenfunction $u$:
$$
-\Delta (u^2 - |du|^2) = 2|Ddu - ug|^2,
\tag 0.2
$$
observed by Lee in [L], 
one may know the scalar
curvature for the conformal compactification $(X, u^{-2}g)$ if one
knows the asymptotic behavior of $u$ near the boundary. This turns
out to be a very interesting construction for its own sake.

The paper is organized as follows. In Section 1 we
will introduce notations and do some computations for the hyperbolic 
spaces. 
In Section 2, we will introduce conformally compact Einstein manifolds 
and relevant properties. And we will apply theory of 
uniformly degenerate elliptic linear PDE to solve for the eigenfunctions 
and their expansions. Finally in Section 3 we will introduce some 
conformal compactifications and prove Theorem 0.1.

\vskip 0.1in
\noindent
{\bf Acknowledgment} \quad 
The author would like to thank Xiao Zhang for bringing my attention 
to the paper [AD] and many stimulating discussions. The author would
also like to thank the referee for many suggestions and advices.

\head 1. On Hyperbolic Spaces \endhead

In this section let us recall some computations on hyperbolic spaces.
We will present in the way that makes our choices of compactification in the
Section 3 become natural. Meanwhile we will establish the notations 
for this paper in this section.
First let us introduce the hyperbolic space as a submanifold in the Minkowski 
space-time. Namely,
$$
H^{n+1} = \{ (x,t) \in R^{n+1,1}: |x|^2 - t^2 = -1, t>0\},
$$
and the hyperbolic metric $g_H$ is induced from the standard Lorenz 
metric $g_L$, i.e.
$$
(H^{n+1}, g_H) = (R^{n+1}, 
\frac{(d|x|)^2}{1+|x|^2} + |x|^2 h_0),
\tag 1.1
$$
where $h_0$ is the standard metric on the unit sphere. It is easily computed
that the coordinate functions in $R^{n+1, 1}$ are all eigenfunctions on 
$H^{n+1}$ for $g_H$, just as the coordinate functions in Euclidean
space are all eigenfunctions on sphere for the round metric. Namely,
$$
\Delta x_i = (n+1)x_i, \ \  \Delta t = (n+1) t, \ \ 
\text{for $i=1, \cdots, n+1$}
\tag 1.2
$$
where $\Delta$ is the trace of the Hessian for the metric $g_H$. 
One conformal compactification near the boundary is 
$$
g_H = \rho^{-2}( (dr)^2 + h_0),
\tag 1.3
$$
where 
$$
\rho = 1/|x|, \ \ \text{and} \ \ \rho = \sinh r.
\tag 1.4
$$
Another way to introduce hyperbolic space is the Poincar\'{e} ball
$$
(H^{n+1}, g_H) = (B^{n+1}, (\frac 2{1-|y|^2})^2|dy|^2),
$$
where $y\in B^{n+1}$ and $B^{n+1}$ is the unit ball in $R^{n+1}$. 
We find that
$$
t = \frac {1+|y|^2}{1-|y|^2},
\tag 1.5
$$
and that $(B^{n+1}, (t+1)^{-2}g_H)$ is the Euclidean ball 
and $(B^{n+1}, t^{-2}g_H)$ is the round hemisphere. 
Finally let us introduce the hyperbolic space by 
the upper half space model
$$
(H^{n+1}, g_H) = (R^{n+1}_+, \frac {|dz|^2}{z_{n+1}^2}),
\tag 1.6
$$
where $z \in R^{n+1}_+$ and $R^{n+1}_+=\{(z_1, z_2, \cdots. z_{n+1})
\in R^{n+1}: z_{n+1} > 0\}$. 
We find that
$$
\Delta \frac 1{z_{n+1}}= (n+1) \frac 1{z_{n+1}},
$$
and $(R^{n+1}_+, (\frac 1{z_{n+1}})^{-2}g_H)$ is the flat half space.
which can be considered as a partial compactification of the hyperbolic
space.

Thus we observe here that the three compactifications: ball, hemisphere 
and half space (partial compactification),
are all through conformal factors which come from 
eigenfunctions. Those are clearly better
than the compactification (1.3) when we need to work not only 
near the boundary.
Before we end this section let us recall the coordinate changes for
the three different models for the hyperbolic space (cf. Chapter 4 in
[R]).
$$
x = \frac 2{1-|y|^2} y, \ \ \ 
\left\{\aligned y_\alpha & = \frac {2z_\alpha}{|z+e_{n+1}|^2} \\
                y_{n+1}  & = \frac {|z|^2 - 1}{|z+e_{n+1}|^2},
       \endaligned\right.
\ \ \ \left\{\aligned z_\alpha & = \frac {2y_\alpha}{|y-e_{n+1}|^2} \\
		z_{n+1}  & = \frac {1 - |y|^2}{|y-e_{n+1}|^2} ,
\endaligned\right.
\tag 1.7
$$
where $e_{n+1} = (0, \cdots, 0, 1)$, and derive
$$
\frac 1{z_{n+1}} = \frac {|y-e_{n+1}|^2}{1 - |y|^2} = 
t - x_{n+1}
\tag 1.8
$$
and
$$
\rho = \frac {1-|y|^2}{2|y|}.
\tag 1.9
$$

\vskip 0.1in
\head 2. Conformally compact Einstein manifolds and the 
positive eigenfunctions. \endhead

Let us first define what is a conformally compact Einstein manifold.
Suppose $X^{n+1}$ is a compact (n+1)-manifold with boundary $\partial X = M^n$.
A Riemannian metric $g$ in the interior of $X$ is said to be 
$C^{m, \alpha}$ conformally compact if $\bar g = r^2 g$ extends as a 
$C^{m, \alpha}$ metric on $\bar X$, where $r$ is a defining function of the
boundary in the sense that: $r>0$ in $X$, $r=0$ and $dr \neq 0$ on $M$. 
Clearly defining functions are not unique. For a given defining function,
the metric $\bar g$ restricted to $TM$ induces a metric $\hat g$ on $M$.
$\hat g$ rescales upon changing the defining function $r$, therefore
defines a conformal class $[\hat g]$ on $M$. We call $(M, [\hat g])$
the conformal infinity of the conformally compact manifold $(X, g)$. 
Conformally compact Einstein manifold $(X, g)$ is a conformally compact
manifold which is Einstein, i.e. $\text{Ric}(g) = -ng$. 

The boundary regularity of the conformally compact Einstein
metric is an important issue. Thanks to M. Anderson [An], in 4 dimension,
we know that $C^{2, \alpha}$ would imply in Theorem 0.1 the full
smoothness of the conformally compact Einstein metric. In other words,
the conformally compact Einstein 4-manifold in Theorem 0.1 is automatically
even analytic. Sufficient boundary regularity is assumed for our
results in higher dimension to ensure the expansion in the following
Lemma 2.1.

Solving a first order nonlinear PDE by the method of characteristics 
introduced in [FG] [GL], one at least has, near the boundary, as follows
(cf. Lemma 5.4 in [AD]). 

\proclaim{Lemma 2.1} 
Given a conformally compact Einstein
manifold with the round sphere as its conformal infinity. Then,
taking the standard round metric $h_0$ and its associated special
defining function, one has
$$
g = \rho^{-2}((dr)^2 + g_r)
$$
where 
$$
g_r = h_0 + \rho^n h, \ \ \text{Tr}_{h_0}h = O(\rho^n),
\ \ \rho = \sinh r.
\tag 2.1
$$
\endproclaim
The choice of the defining function $\rho$ in this lemma
is different from the
choice made by Fefferman and Graham in [FG], particularly in 
Lemma 2.1 in [Gr]. Because of this choice the expansion of
$g_r$ becomes much nicer. And this choice comes from (1.3).
One thing we learn from Section 1 is that the positive eigenfunctions 
seem to play a role in conformal compactifications. 
Therefore, let us turn our attention to
study the positive eigenfunctions. We will 
use the theory of uniformly degenerate
elliptic linear PDE on conformally compact manifolds developed 
in [M], [L] and [GL]. We first recall 
an analytic lemma with modifications for simplicity in the
following Lemma 2.2 from [M]. [L] and [GL] 
(cf. Proposition 3.3 in [L] for example). 
Given a conformally compact manifold with a fixed defining function $r$ 
which may be defined globally and is identical to $r$ as in Lemma 2.1 
near the boundary. We may define the space of functions 
$$
C^{m, \alpha}_s(X) = \rho^s C^{m,\alpha}(X, g)
\tag 2.2
$$
where $C^{m,\alpha}(X,g)$ is defined as usual for a Riemannian manifold and
$\alpha \in (0,1)$. And 
$$
\|f\|_{C^{m, \alpha}_s(X)} = \|\rho^{-s}f\|_{C^{m,\alpha}(X, g)}.
\tag 2.3
$$
In our situation when the conformally compact manifold is Einstein with
the round sphere as its conformal infinity, we have a coordinate near the
boundary as the hyperbolic space, an annular region in $y$ coordinates
(cf. please see notations in Section 1), for example.

\proclaim{Lemma 2.2} 
Suppose that $(X, g)$ is a conformally compact manifold. Then 
$$
- \Delta + (n+1): C^{m+2, \alpha}_s \longrightarrow C^{m,\alpha}_s
\tag 2.4
$$
is an isomorphism whenever
$$
-1 < s < n+1
\tag 2.5
$$
\endproclaim

Next we are going to find eigenfunctions on 
the conformally compact Einstein manifold with the round sphere as its
conformal infinity. It makes sense to consider functions
$t, x_i$ for $i = 1, \cdots, n+1$ near the boundary because we may
adopt for $X$ near the boundary the same coordinate systems 
that we had in Section 1 for the hyperbolic space. 
To simplify the notation we will use $f = O(\rho^s)$ to stand for
$f \in C^{2, \alpha}_s (X)$.

\proclaim{Lemma 2.3} Suppose that $(X, g)$ is a conformally compact
Einstein manifold with the round sphere as its conformal infinity. 
Then there exist eigenfunctions $u>0$, and $v_1, v_2, \dots, v_{n+1}$
such that
$$
\Delta u = (n+1)u, \ \Delta v_i = (n+1) v_i, \ \text{for} \ 
i = 1, 2, \dots, n+1
\tag 2.6
$$
and
$$
u = t + O(\rho^{n+\mu}) \ \text{and} \ 
v_i = x_i + O(\rho^{n+\mu}) \ \text{for any} \ \mu < 1.
\tag 2.7
$$
\endproclaim
\demo{Proof}  First let us construct $u$ by the above
Lemma 2.2. Since we are using the same coordinate systems as we
did for the hyperbolic space, the only thing that is different from
the hyperbolic space is the metric $g_r$. We calculate
$$
\Delta t = \frac {\rho^{n+1}}{\sqrt{\det g_r}}\partial_r(
\rho^{1-n} \sqrt{\det g_r}\partial_r t)
= (n+1)t - \frac 12 \text{Tr}_{h_0}g_r'
= (n+1)t + O(\rho^{2n-1}),
$$
where $g_r = h_0 + \rho^n h$ according to Lemma 2.1, 
the prime means differentiating with respect to $r$,
and note also that $t' = -1/\rho^2$.
Then, by Lemma 2.2, we know there exists a positive
function $u$ which satisfies $\Delta u = (n+1)u$ and
$$
u = t + O(\rho^{n+\mu}) \ \text{for any} \ \mu < 1.
$$
Similarly, let us compute
$$
\Delta x_i = \frac {\rho^{n+1}}{\sqrt{\det g_r}}\partial_r(
\rho^{1-n} \sqrt{\det g_r} \partial_r x_i) + 
\rho^2 \frac 1{\sqrt{\det g_r}}\partial_\gamma (
\sqrt{\det g_r}g^{\gamma\delta}_r\partial_\delta x_i),
$$
where 
$$
\frac {\rho^{n+1}}{\sqrt{\det g_r}}\partial_r(
\rho^{1-n} \sqrt{\det g_r} \partial_r x_i)
= \frac {\rho^{n+1}}{\sqrt{\det h_0}}\partial_r(
\rho^{1-n} \sqrt{\det h_0}\partial_r x_i) + \frac 12 \rho^2 x_i'
\text{Tr}_{h_0}g_r'
$$
and
$$
\aligned
& \frac {\rho^2}{\sqrt{\det g_r}} \partial_\gamma (
\sqrt{\det g_r}g^{\gamma\delta}_r\partial_\delta x_i) \\
& = \frac{\rho^2}{\sqrt{\det g_r}}\partial_\gamma (
\sqrt{\det g_r}h^{\gamma\delta}_0\partial_\delta x_i) 
+ \frac {\rho^{n+1}}{\sqrt{\det g_r}}\partial_\gamma (
\sqrt{\det g_r}h^{\gamma\delta}\partial_\delta (\rho x_i)) \\
& =  \frac {\rho^2}{\sqrt{\det h_0}}\partial_\gamma (
\sqrt{\det h_0}h^{\gamma\delta}_0\partial_\delta x_i).
+ \frac 12 \rho^2 h_0^{\gamma\delta}\text{Tr}_{h_0}
\partial_\delta (g_r - h_0) \\
& \ \ \ \ \ + \rho^{n+1}\frac 1{\sqrt{\det g_r}}\partial_\gamma (
\sqrt{\det g_r}h^{\gamma\delta}\partial_\delta (\rho x_i)).
\endaligned
$$
Therefore
$$
\Delta x_i = 
=  (n+1) x_i + O(\rho^{n+1}).
\tag 2.8
$$
Therefore, applying Lemma 2.2 again,
we obtain functions $v_i$ which solves $\Delta v_i = (n+1)v_i$ and
$$
v_i = x_i + O(\rho^{n+\mu}), \ \ \text{for any} \ \mu < 1.
$$
\enddemo

\vskip 0.1in
\head 3. Conformal compactifications. \endhead

In this section we will
conformally compactify the manifolds with eigenfunctions obtained
in Lemma 2.3. 
Before going to the proof of Theorem 0.1, let us mention that 
the compactification 
$(X, u^{-2}g)$ obtained by the eigenfunction $u$ in Lemma 2.3 
is a compact manifold with a totally geodesic standard
sphere boundary and scalar curvature $\geq n(n+1)$.
This can be easily verified similar to what we will do in 
the proof of following Lemma 3.2. It is worthwhile to note that the 
Bochner formula 
$$
-\Delta (u^2-|\nabla u|^2) = 2|Ddu - ug|^2
\tag 3.1
$$
for the eigenfunction $u$ may very well 
invite people to show that $u^2- |du|^2 = 1$ and $Ddu = ug$, which 
would quickly imply the rigidity. But, apparently, estimate (2.7) from 
Lemma 2.2 miserably just fails
to provide the sufficient decay of $u^2 - |du|^2 -1$. 

To prove our Theorem 0.1 
we consider the partial compactification corresponding to the half 
space. 
Following the notations in Section 1, we are changing into
upper half space, $z$ coordinate.
The metric in this coordinate becomes
$$
g = \frac 1{z_{n+1}^2}|dz|^2 + \rho^{n-2} h
$$
in the light of (1.3) and (2.1). 
We first pay attention to the tail term $\rho^{n-2} h$ in $z$
coordinate.

\proclaim{Lemma 3.1} In $z$ coordinates, 
$$
\rho^{n-2} h = |z|^{-n-2} \tilde h_{ij}dz_i dz_j
\tag 3.2
$$
where $\tilde h_{ij}$ are well bounded in the sense that
$$
|\tilde h_{ij}| + |z||\partial_z \tilde h_{ij}| + |z|^2 
|\partial^2_z \tilde h_{ij}| < \infty, 
\tag 3.3
$$
at least for $|z|$ very large.
\endproclaim
\demo{Proof} The proof is simply to perform the coordinate
change. Again we follow notations used in Section 1.
First, it is easily seen that $|y-e_{n+1}|$ is very small
when $|z|$ is very large and vice verse. So let us restrict 
ourselves to
the very small neighborhood of $e_{n+1}$. Then, according to (1.7),
$$
\rho = \frac 1{|x|} = \frac {1-|y|^2}{2|y|} = 
\frac {z_{n+1}|y-e_{n+1}|^2}{2|y|}
= \frac {2z_{n+1}}{|z+e_{n+1}|^2|y|}.
\tag 3.4
$$
Meanwhile, in $y$ coordinates,
$$
h = h_{\alpha\beta}d\phi_\alpha d\phi_\beta,
$$
where $\{ d\phi_\alpha \}_{\alpha=1}^{n}$ is an orthonormal co-frame on the 
unit sphere with the round metric
and $h_{\alpha\beta} \in C^2$. We may write
$$
d\phi_\alpha = \sum_{i=1}^{n+1}c^\alpha_i dy_i,
\tag 3.5
$$
where $\{c^\alpha_i \}_{i=1}^{n+1}$ are all well bounded.
Therefore, by the transformation formula (1.7)
$$
\left\{\aligned dy_\alpha & = \frac 2{|z+e_{n+1}|^2}dz_\alpha - 
\frac 4{|z+e_{n+1}|^4}z_\alpha (z+e_{n+1})dz \\
dy_{n+1} & =  \frac 2{|z+e_{n+1}|^2}  
\frac {2z_{n+1}+2}{|z+e_{n+1}|^2} zdz 
- \frac {2(|z|^2 -1)}{|z+e_{n+1}|^4} dz_{n+1}
\endaligned\right.
$$
This implies
$$
\frac {\partial y_i }{\partial z_j} = O(|z|^{-2}).
\tag 3.6
$$
Thus
$$
\rho^{n-2}h =  |z|^{-n-2} \tilde h_{ij} dz_i dz_j,
$$
where $\tilde h_{ij}$ are well bounded as desired. So the lemma
is proved.
\enddemo

We now consider the positive eigenfunction
$$
\psi = u - v_{n+1}
= t - x_{n+1} + O(\rho^{n+\mu})=
\frac 1{z_{n+1}} + O(\rho^{n+\mu}) \ \ \text{for any $\mu <1$} 
\tag 3.7
$$
in the light of (1.8) in Section 1.
Again, positivity comes from a maximum principle and the 
boundary behavior of $u - v_{n+1}$. More importantly 
we have
\proclaim{Lemma 3.2} The scalar
curvature of the new metric $g_h = \psi^{-2}g$ 
is nonnegative and integrable.
\endproclaim

\demo{Proof} First we calculate the scalar curvature for 
the metric $\psi^{-2} g$
$$
\Delta \psi^{-\frac {n-1}2}
= -\frac {n-1}2 \psi^{-\frac {n+1}2}
\Delta \psi + \frac {n^2-1}4 \psi^{-\frac {n+3}2} |\nabla \psi|^2,
$$
that is
$$
- \Delta \psi^{-\frac {n-1}2} - \frac {n^2-1}4 \psi^{-\frac {n-1}2}
= \frac {n^2-1}4 (\psi^2 - |\nabla\psi|^2)\psi^{-\frac {n+3}2}.
\tag 3.8
$$
Therefore
$$
R[\psi^{-2}g] = n(n+1)(\psi^2 - |\nabla\psi|^2).
$$
Recall the Bochner formula for the eigenfunctions observed in [L]
for $\psi$
$$
-\Delta (\psi^2 - |\nabla\psi|^2) = 2|Dd\psi - \psi g|^2.
\tag 3.9
$$
Thus to prove the scalar curvature $R[\psi^{-2}g]\geq 0$ one only
needs to apply a maximum principle and to verify that 
$\psi^2 - |\nabla\psi|^2$ goes to zero towards the infinity. In fact
we have
$$
\aligned
(u-v_{n+1})^2 & - |du - dv_{n+1}|^2
= O(\frac {\rho^{n+\mu}}{z_{n+1}}) \\
= & O(|y-e_{n+1}|^2 \rho^{n-1+\mu}) = O (|z|^{-n-1-\mu})
\endaligned
\tag 3.10
$$
in the light of (1.7), (2.7), (3.4) and (3.7). (3.10) also implies
that the scalar curvature is integrable with respect to the metric $g_h$.
It turns out that 
(3.10) is another key calculation in our approach regarding the
remark made right after the formula (3.1).
\enddemo

For the convenience of readers we recall the definition of an asymptotically
flat manifold. We simply use Definition 6.3 in [LP].

\proclaim{Definition 3.3} (Definition 6.3 [LP])
A Riemannian $(n+1)$-manifold $(M^{n+1}, g)$ is an asymptotically flat
manifold of order $\tau$ if there exists a decomposition $M = M_0 \bigcup
M_\infty$ (with $M_0$ compact) and a diffeomorphism 
$M_\infty \longleftrightarrow
R^{n+1}\setminus B_R$ for some $R>0$, satisfying, if $(M, g) = 
(R^{n+1}\setminus B_R, g_{ij})$,
$$
g_{ij} = \delta_{ij} + O(|z|^{-\tau}), \ 
\partial_k g_{ij} = O(|z|^{-\tau-1}), \ 
\partial_k\partial_l g_{ij} = O(|z|^{-\tau-2})
$$
for all $i,j,k,l = 1, \dots, n+1$, 
as $|z|\longrightarrow \infty$ in the asymptotic coordinate chart 
$R^{n+1}\setminus B_R$.
\endproclaim

Also, given an asymptotically flat $(n+1)$-manifold and an 
asymptotic coordinate
chart, one may define the mass for the asymptotically 
flat $(n+1)$-manifold as
follows: (cf. Definition 8.2 in [LP])
$$
m(g) = \lim_{R \longrightarrow\infty} \frac 1{|S^n|}
\int_{S_R} \sum_{i,j=1}^{n+1}
(\partial_ig_{ij}-\partial_jg_{ii})\frac{z_j}{|z|} d\sigma
$$
where $|S^n|$ is the volume of the unit n-sphere, if this limit exists.

\proclaim{Lemma 3.4} We may consider the doubling $(Y, G)$ of $(X, g_h)$ 
along the boundary in above partial compactification. Then $(Y, G)$
is at least $C^{n-1,1}$ (more precisely, $C^{n,1}$ if 
$n$ is even and $C^{n-1,1}$ is $n$ is odd), and is an
asymptotically flat manifold of order $n-1+\mu$, therefore $m(G)=0$, 
with integrable nonnegative scalar curvature.
\endproclaim

\demo{Proof} By the above Lemma 3.2 one knows that
the scalar curvature for $G$ is nonnegative. Note $G$ is 
$C^{2,1}$ at least, therefore its curvature tensor is well-defined.

To understand the metric $G$ better, we have
$$
g_h = \psi^{-2}g = \psi^{-2}\frac 1{z_{n+1}}|dz|^2 
+ \psi^{-2}\rho^{n-2}h.
\tag 3.11
$$
Let us first recall
$$
\psi z_{n+1} = 1 + z_{n+1}O(\rho^{n+\mu}),
$$
then, plugging into (3.11), 
$$
\aligned
g_h  = & \frac 1{(1+ z_{n+1}O(\rho^{n+\mu}))^2} |dz|^2 
+ \psi^{-2}\rho^{n-2}h \\
= & |dz|^2 + z_{n+1}O(\rho^{n+\mu})|dz|^2 + z_{n+1}^2 \rho^{n-2}
\frac 1{|z+e_{n+1}|^4}\tilde h_{ij}dz_i dz_j. 
\endaligned
$$
When $(z_1, \dots, z_n)$ is fixed, from (3.4),
$$
g_h = |dz|^2 + O(z_{n+1}^{n+1+\mu})|dz|^2 + O(z_{n+1}^n)\tilde h_{ij}
dz_i dz_j
\tag 3.12
$$
near the boundary, which tells us the smoothness of the metric $G$
at the 
boundary $z_{n+1} =0$. 
Meanwhile, when $|z|$ is very large, we have, from (3.4),
$$
g_h = |dz|^2 + O (|z|^{-(n-1+\mu)})|dz|^2 + O(|z|^{-n})\tilde h_{ij}
dz_i dz_j.
\tag 3.13
$$

To check if $(Y, G)$ is an asymptotically flat metric one needs
to verify that all terms $O(|z|^{-k})$ in (3.13) in fact satisfy
$\partial_z O(|z|^{-k}) = O(|z|^{-k-1})$ and $\partial^2_z O(|z|^{-k})
= O(|z|^{-k-2})$. Those indeed are true according to Lemma 2.2 for
$m\geq 2$, Lemma 3.1, and (3.6). For example,
$$
\partial_z f = \frac {\partial y}{\partial z} \partial_y f,
$$
therefore
$$
\partial_z f = O(|z|^{-2})\partial_y f
= O(|z|^{-2}) O(\rho^{k-1}) = O(|z|^{-k-1})
$$
if $f = O(\rho^k)$. Thus the lemma is proved.
\enddemo

\proclaim{Theorem 3.5} Suppose that $(X^{n+1},g)$ is a conformally 
compact Einstein manifold
with the round sphere as its conformal infinity, and
$3\leq n\leq 6$. Then $(X^{n+1},g)$ must be 
the hyperbolic space.
\endproclaim
\demo{Proof} This is a more or less 
straight consequence of the positive mass theorem of
Schoen and Yau [Sc]
except we need to make sure that their theorem applies to 
less smooth metrics as ours. 
By Theorem 4.2 in [Sc], for example, $(Y, G)$ has to 
be isometric 
to $R^{n+1}$, since $(Y, G)$ is an asymptotically flat manifold 
of order $n-1+\mu$ with integrable nonnegative scalar curvature
because of (3.10), and zero mass by Lemma 3.4. In the following we will
point out that the positive mass 
theorem of Schoen and Yau indeed works for asymptotic manifolds which
are at least $C^{2,1}$. First it is easily seen that
Proposition 4.1 in [Sc] still holds in our situations without much
modifications. 
Since the metric $\bar g$ in the proof of Proposition 4.1 in [Sc] was
constructed from cutoff metric $g^{(\sigma)}$ which is Euclidean near
the infinity, one may assume
the asymptotically flat metrics are $C^\infty$ near the infinity
to show that the mass have to be nonnegative following the minimal 
hypersurface argument in [Sc] (dimension is assumed to be less than 
and equal to seven in [Sc]). Then to show that the mass
is zero implies that the asymptotically flat manifold has to be Euclidean
one may follow the proof of 
Lemma 10.7 in [LP] (see also Lemma 3 and Proposition 3 
in [Sc1]). The key is to show that the mass of $G$ is zero implies that 
$G$ is Ricci flat, in the light of Proposition 10.2 in [LP]. 
But in the argument given on page 84-85 in [LP] for
Lemma 10.7 and in the argument given on page 80-81 in [LP] for 
Proposition 10.2 only derivatives of the metrics up to the
second order were involved. 
And the minimality of the zero mass still holds since
we just showed that the mass has to be nonnegative, for dimension
less than 8. 
Moreover the variational formula (8.11) in [LP] certainly holds
for metrics of $C^{2,1}$. Thus with little modifications 
Lemma 10.7 in [LP] holds in our cases.

\enddemo

\proclaim{Remark 3.6}
The application of positive mass theorem to the doubling manifolds 
has been used by Escobar in [Es] in his proof of the Yamabe 
problem for manifolds with boundary. In the 
appendix of the paper [Es] he explained how one can apply 
the positive mass theorem to the doubling manifold which is 
asymptotically flat. 
\endproclaim

\bigskip 
\noindent
{\bf References}:
\roster

\vskip 0.1in
\item"{[AD]}" L. Andersson and M. Dahl, Scalar curvature rigidity for
asymptotically hyperbolic manifolds, Ann. Global Anal. and Geo. 16(1998)
1-27.

\vskip 0.1in
\item"{[An]}" M. Anderson, Boundary regularity, uniqueness and non-uniqueness 
for AH Einstein metrics on 4-manifolds, Preprint. 

\vskip 0.1in
\item"{[An1]}" M. Anderson, Einstein metrics with prescribed conformal
infinity on 4-manifolds, Preprint.


\vskip 0.1in
\item"{[CH]}" Piotr Chru\'{s}ciel and M. Herzlic, The mass of
asymptotically hyperbolic Riemannian manifolds, Preprint
math.DG/0110035.

\vskip 0.1in
\item"{[E]}" J. Escobar, The Yamabe problem on manifolds with boundary,
J. Diff. Geo. 35(1992) 21-84.

\vskip 0.1in
\item"{[FG]}" C. Fefferman and C. R. Graham, Conformal invariants, in
Elie Cartan et les Mathematiques d'Aujourd'hui, Asterisque (1985)
95-116.

\vskip 0.1in
\item"{[Gr]}" C. R. Graham, Volume and Area renormalizations for
conformally compact Einstein metrics. The Proceedings of the 19th Winter School
"Geometry and Physics" (Srn\`{i}, 1999).
Rend. Circ. Mat. Palermo (2) Suppl. No. 63 (2000), 31--42.

\vskip 0.1in
\item"{[GL]}" C.R. Graham and J. Lee, Einstein metrics with prescribed conformal
infinity on the ball. Adv. Math. 87 (1991), no. 2, 186--225.

\vskip 0.1in
\item"{[L]}" John Lee, 
The spectrum of an asymptotically hyperbolic Einstein manifold. Comm.
Anal. Geom. 3 (1995), no. 1-2, 253--271.

\vskip 0.1in
\item"{[LP]}" J, Lee and T. Parker, The Yamabe problem, Bull. Amer. Math. Soc.
17(1987) 37-91.

\vskip 0.1in
\item"{[M]}" R. Mazzeo. The Hodge cohomology of a conformally compact 
metric, J. Diff. Geom. 28(1988) 309-339.

\vskip 0.1in
\item"{[R]}" John G. Ratcliffe, ``Foundations of hyperbolic manifolds'',
Graduate Texts in Mathematics 149, Springer-Verlag, New York, 1994.

\vskip 0.1in
\item"{[Sc]}" R. Schoen, Variational theory for the total scalar curvature
functional for Riemannian metrics and related topics. 
Topics in calculus of variations (Montecatini Terme, 1987),
120--154, Lecture Notes in Math., 1365, 
\newline Springer, Berlin, 1989.

\vskip 0.1in
\item"{[Sc1]}" R. Schoen, Conformal deformation of a Riemannian 
metric to a constant scalar curvature, J. Diff. Geo. 20 (1984)
479-495.

\vskip 0.1in
\item"{[Wa]}" X. Wang, Uniqueness of AdS space-time in any dimension,
Preprint Spet. 2002.

\vskip 0.1in
\item"{[W]}" E. Witten, A new proof of the positive energy theorem, Comm.
Math. Phys. 80(1981) 381-402.

\vskip 0.1in
\item"{[Z]}" X. Zhang, A definition of total energy-momenta and the
positive mass theorem for asymptotically hyperbolic 3-manifolds, Preprint
2001.

\endroster
\enddocument